\documentclass[a4paper,12pt]{article}

\newtheorem{theorem}{Theorem}

\usepackage{amsmath}
\usepackage{amssymb}
\usepackage{amsfonts}
\usepackage{mathrsfs}
\usepackage{graphicx}

\renewcommand{\b}{\mathbf{b}}
\renewcommand{\c}{\mathbf{c}}
\newcommand{\p}{\mathbf{p}}
\renewcommand{\v}{\mathbf{v}}
\newcommand{\y}{\mathbf{y}}
\newcommand{\x}{\mathbf{x}}
\renewcommand{\y}{\mathbf{y}}
\newcommand{\z}{\mathbf{z}}

\newcommand{\A}{\mathbf{A}}
\newcommand{\M}{\mathbf{M}}
\newcommand{\X}{\mathbf{X}}
\newcommand{\Y}{\mathbf{Y}}

\newcommand{\R}{\mathbb{R}}
\renewcommand{\S}{\mathbb{S}}
\newcommand{\Mom}{\mathscr{M}}
\newcommand{\Pos}{\mathscr{P}}
\newcommand{\Spt}{\mathscr{X}}

\usepackage{color}
\usepackage{ulem}
\newcommand{\new}[1]{{\color{black}#1}}

\title{Dual optimal design and the Christoffel-Darboux polynomial}

\begin{document}
	
	\author{Yohann De Castro$^1$, Fabrice Gamboa$^2$, Didier Henrion$^{3,4}$, Jean Bernard Lasserre$^{2,3}$}
	
	\footnotetext[1]{Institut Camille Jordan UMR 5208, \'Ecole Centrale de Lyon, 36 Avenue Guy de Collongue, F-69134 Écully, France.}
	\footnotetext[2]{Institut de Math\'ematiques de Toulouse, Universit\'e Paul Sabatier, CNRS, 118 route de Narbonne, F-31062 Toulouse, France.}
	\footnotetext[3]{CNRS, LAAS, 7 avenue du colonel Roche, F-31400 Toulouse, France.}
	\footnotetext[4]{Faculty of Electrical Engineering, Czech Technical University in Prague,
		Technick\'a 2, CZ-16626 Prague, Czechia.}
	
	\date{Draft of \today}
	
	\maketitle
	
\begin{abstract}
The purpose of this short note is to show that the Christoffel-Darboux polynomial, useful in approximation theory and data science, arises naturally when deriving the dual to the problem of semi-algebraic D-optimal experimental design in statistics. It uses only elementary notions of convex analysis. \new{Geometric interpretations and algorithmic consequences are mentioned.}
\end{abstract}

\section{Introduction}

In \cite{aos19} the problem of optimal design of statistical experiments was revisited in the broad framework of polynomial regressions on semi-algebraic domains. A numerical solution was proposed, based on the so-called moment-SOS (sums of squares) hierarchy of semidefinite programming relaxations \cite{l10}. While optimality arguments were used in \cite{aos19} to derive many of the results, the dual to the problem of optimal experimental design was not explicitly constructed and studied. It is the purpose of this note to clarify this point in a self-contained and direct way. \new{We believe that a significantly shorter, separate, elementary derivation of the dual is welcome, and its interpretation is informative.} 
	
We use elementary arguments of convex analysis to show how the Christoffel-Darboux polynomial, ubiquitous in approximation theory and data science \cite{lp19}, arises naturally in the dual to the D-optimal design problem.

\section{\new{Primal formulation of the approximate design problem}}

Let  $\S^n$ denote the set of symmetric real matrices of size $n$. Given a vector $\v_d(\x)$ whose elements form a basis of the vector space of real polynomials of degree up to $d$ in the vector indeterminate $\x \in \R^n$, let
\[
\Mom_d(\Spt):=\left\{\y \in \R^{n_d} : \y = \int_\Spt \v_d(\x)d\mu(\x) \:\text{for some positive Borel measure $\mu$ on $\Spt$}\right\}
\]
denote the convex cone of moments of degree up to $d$ on a given compact semi-algebraic set
$\Spt \in \R^n$ with non-empty interior. This cone has dimension
\[
n_d:=\frac{(n+d)!}{n!\:d!}.
\]
Given a measure $\mu$ and its moment vector $\y \in \Mom_{2d}(\Spt)$, define the moment matrix
\[
\M_d(\y):=\int_\Spt \v_d(\x)\v^\star_d(\x)d\mu(\x) \in \S^{n_d}
\]
where the star denotes transposition. Each entry of the above matrix $\v_d(\x)\v^\star_d(\x)$ is a polynomial of degree up to $2d$, and hence it is a linear combination of elements in basis vector $\v_{2d}(\x)$. Consequently, each entry of $\M_d(\y)$ is a linear combination of entries of the moment vector $\y$, whose dimension is $n_{2d}$. It follows that $\M_d$ can be interpreted as a linear map from $\R^{n_{2d}}$ to $\S^{n_d}$, 
and let $\M^\star_d$ denote its adjoint map from $\S^{n_d}$ to $\R^{n_{2d}}$, defined such that
\begin{equation}\label{adjoint}
\text{trace}(\M_d(\y)\X) = \M^\star_d(\X)\y
\end{equation}
holds for all $\X \in \S^{n_d}$ and $\y \in \R^{n_{2d}}$.

Given a matrix $\A \in \R^{m\times n_d}$, a vector $\b \in \R^m$, a vector $\c \in \R^{n_d}$ and a strictly concave function $\phi$ from $\S^{n_d}$ to $\R$, consider the primal optimal design problem
\begin{equation}\label{primal}
\begin{array}{ll}
\inf_{\y} & \c^\star  \y - \phi(\M_d(\y))\\
\text{s.t.} & \A\y = \b \\
& \y \in \Mom_{\new{2}d}(\Spt)
\end{array}
\end{equation}
where the infimum is with respect to vectors $\y \in \R^{n_d}$.

\section{\new{Dual formulation of the approximate design problem}}

In problem (\ref{primal}) we introduce the matrix variable $\Y \in \S^{n_d}$ and equality constraint $\Y=\M_d(\y)$, the Lagrange multipliers $\X$, $\z$, $\p$ and the Lagrangian
\[
f(\X,\z,\p,\Y,\y):=\c^\star \y - \phi(\Y) + \text{trace}\left\{\X(\Y-\M_d(\y))\right\} + \z^\star (\b-\A\y) - \p^\star \y
\]
where $\X \in \S^{n_d}$, $\z \in \R^m$ and $\p$ belongs to $\Pos_{\new{2}d}(\Spt)$, the convex cone of polynomials \new{of degree at most $2d$} that are positive on $\X$, which is dual to the moment cone $\Mom_{\new{2}d}(\Spt)$ according to the Riesz-Haviland Theorem \cite[Theorem 3.1]{l10}.
Use \eqref{adjoint} to rearrange terms as follows
\[
f(\X,\z,\p,\Y,\y) = \b^\star  \z + \left\{\text{trace}(\X\Y)-\phi(\Y)\right\} +  \{\c^\star  - \z^\star \A - \M^\star_d(\X) -\p^\star\}\y.
\]
The dual problem to \eqref{primal} is obtained by minimizing the dual function
\[
g(\X,\z,\p):=\inf_{\Y\in \S^{n_d},\:\y\in\R^{n_{\new{2}d}}} f(\X,\z,\p,\Y,\y).
\]
First observe that this function is bounded below only if
\begin{equation}\label{poly}
\p^\star = \c^\star  - \z^\star \A - \M^\star_d(\X).
\end{equation}
Defining
\[
\phi^\star(\X) := \inf_{\Y\in\S^{n_d}} \left\{\text{trace}(\X\Y)-\phi(\Y)\right\}
\]
as the concave conjugate function to $\phi$, the problem of maximizing the dual function becomes
\begin{equation}\label{dual}
\begin{array}{ll}
\sup_{\X,\z} & \b^\star \z + \phi^\star(\X) \\
\text{s.t.} & \c^\star - \z^\star \A - \M^\star_d(\X) \in \Pos_{\new{2}d}(\Spt)
\end{array}
\end{equation}
where the supremum is with respect to matrices $\X \in \S^{n_d}$ and vector $\z \in \R^m$.
\new{Letting $p(x):= \M^\star_d(\X)\v_{2d}(\x) = \v_d(\x)^\star \X \v_d(\x)$},
the conic constraint in dual problem \eqref{dual} can be formulated
as a polynomial positivity constraint
\[
(\c^\star-\z^\star \A)\v_{2d}(\x) \geq p(\x)
\]
satisfied for all $\x \in \Spt$.

\begin{theorem}[Duality for optimal design]\label{duality}
Problems \eqref{primal} and \eqref{dual} are in strong duality, i.e. their values coincide.
\end{theorem}

{\bf Proof:}
Weak duality, i.e. the value of primal problem \eqref{primal} is greater than or equal to the value of dual problem \eqref{dual}, follows from the inequality
\[
\phi^\star(\X)+\phi(\Y) \leq \text{trace}(\X\Y)
\]
which holds by definition of $\phi^\star$ for every positive semidefinite pair $\X$, $\Y$.
Indeed, for any feasible $\X$, $\y$, $\z$ it holds
\[
0 \leq (\c^\star  - \z^\star\A - \M^\star_d(\X))\y = \c^\star\y -\b^\star\z - \text{trace}(\X\M_d(\y))
\leq \c^\star\y  - \b^\star\z - \phi^\star(\X) - \phi(\M_d(\y))
\]
and hence
\[
\c^\star\y - \phi(\M_d(\y)) \geq \b^\star\z + \phi^\star(\X).
\]
Strong duality, i.e. the above inequality is an equality for any optimal values $\hat{\X}$, $\hat{\y}$, $\hat{\z}$ follows from concavity of function $\phi$ and the so-called Slater qualification constraint \cite[Section C.1]{l10}, i.e. the existence of an interior point for primal problem \eqref{primal}: choose e.g. the vector of moments of an atomic measure supported on $\Spt$ with more than $n_d$ distinct atoms. Then the moment matrix $\M_d(\y)$ is positive definite. Equivalently, the complementarity condition
\begin{equation}\label{comp}
\hat{\p}^\star \hat{\y} = 0
\end{equation}
holds for $\hat{\p}^\star:=\c^\star  - \hat{\z}^\star\A - \M^\star_d(\hat{\X})$ as in \eqref{poly}.
$\Box$

\section{Christoffel-Darboux polynomial}

In \cite{aos19} various functions $\phi$ are considered, depending on the optimal design problem of interest.
In optimal design problem \eqref{primal}, let
\begin{equation}\label{logdet}
\phi(\Y) := \log\det\Y, \quad \A \y := \int_\Spt d\mu(\x), \quad \b = 1, \quad \c = 0
\end{equation}
i.e. we are minimizing over probability measures supported on $\Spt$ a function $\phi$ which is
the classical barrier function used in interior point methods for semidefinite programming \cite{nn94}. The domain of $\phi$ is the cone of positive definite matrices.
This optimal design problem has the same solution as the D-optimal design problem corresponding to the positively homogeneous objective function $(\det\Y)^{-1/n_d}$.

\begin{theorem}\label{cd}
Problem \eqref{primal} with data \eqref{logdet} has a unique solution
\[
\hat{\y} \in \Mom_{\new{2}d}(\{\x \in \Spt : p(\x) = n_d\})
\]
where $p(\x):=\v^\star_d(\x)\M^{-1}_d( \hat{\y})\v_d(\x)$ is the Christoffel-Darboux polynomial \new{associated to $\hat{\y}$.}
\end{theorem}

{\bf Proof:}
Uniqueness of the solution $\hat{\y}$ follows from convexity of the feasibiliy set and strict concavity of the objective function in problem \eqref{primal}.
The solution $\hat{\y}$ is an interior point, i.e. $\M_d(\hat{\y})$ is positive definite.  As explained in the proof of Theorem \ref{duality}, the Karush-Kuhn-Tucker (KKT) optimality conditions
are necessary and sufficient for an optimal solution, see e.g. \cite[Section C.1]{l10}: all partial derivatives of the Lagrange dual function $f$ must vanish, and this implies that
\[
\frac{\partial f}{\partial \Y} = \hat{\X}-\frac{\partial \phi}{\partial \Y}(\M_d(\hat{\y})) = 0
\]
and hence that
\[
\hat{\X} = \M^{-1}_d(\hat{\y})
\]
for an optimal primal-dual pair $\hat{\X}$, $\hat{\y}$. From the complementarity condition \eqref{comp} and property \eqref{adjoint} we deduce that the optimal $\hat\z$ satisfies
\[
\hat{\z} = -\M^\star_d(\M^{-1}_d(\hat{\y}))\hat{\y} = -\text{trace}(\M_d(\hat{\y})\M^{-1}_d(\hat{\y})) = -n_d.
\]
Complementarity condition \eqref{comp} means that an optimal vector of moments $\hat{\y}$ corresponds to a measure $\mu$ supported on the zero level set of the optimal
positive polynomial with coefficients $\hat{\p}:=\M^\star_d(\M^{-1}_d(\hat{\y}))-n_d$, i.e. the algebraic set $\{\x \in \Spt : p(\x) = n_d\}$.$\Box$

The concave conjugate function is $\phi^\star(\X)=n_d+\phi(\X)$, its domain is the cone of positive definite matrices, 
and from the proof of Theorem \ref{cd}, the dual design problem (\ref{dual}) has the simple form
\begin{equation}\label{dual-final}
\begin{array}{ll}
\sup_\X &  \log\det\X \\
\text{s.t.} & n_d - \v_d(\x)^\star \X \v_d(\x) \geq 0, \quad \forall \x \in \Spt.
\end{array}
\end{equation}
Its solution is $\hat{\X}=\M^{-1}_d(\hat{\y})$ where $\hat{\y}$ is the unique solution of problem (\ref{primal}).
The Christoffel-Darboux polynomial $p(\x)=\v^\star_d(\x)\hat{\X}\v_d(\x)$ is SOS since matrix $\hat{\X}$ is positive definite, and dual matrix $\X$ is such that $n_d \geq \v_d(\x)^\star \X \v_d(\x) \geq 0$ for all $\x \in \Spt$.

From Theorem \ref{cd} the optimal sequence of moments $\hat{\y}$ in primal problem \eqref{primal} has a representing atomic measure $\mu$ whose atoms are given by  the level set of the Christoffel-Darboux polynomial. This sequence of moments is unique, but there could be another measure, atomic or not, with the same moments.

Dual problem \eqref{dual-final} has also an interpretation in computational geometry. Indeed if $d=1$ then 
$\x\mapsto p(\x)=n_d-\v_d(\x)^T\X\v_d(\x)$ is a quadratic  polynomial and so the set 
$\mathcal{E}:=\{\x \in \R^n : p(\x)\geq0\}$  is an ellipsoid that contains $\Spt$, and $\log\det\X$ is related to the volume of $\mathcal{E}$.
So the dual design problem is also equivalent to the problem of finding the ellipsoid of minimum volume that contains $\Spt$, which is the celebrated L\"owner-John ellipsoid problem. For $d=1$ this was already observed in \cite{t16} and therefore 
\eqref{dual-final} can be considered as a generalization to the case $d>1$ of the L\"owner-John ellipsoid problem
with set $\Spt$ and  with $\log\det \X$ as a proxy for the volume of $\mathcal{E}$.

\section{Conclusion}

In this note we use only elementary concepts of convex analysis to show that the Christoffel-Darboux polynomial, so useful in approximation theory and data analysis  \cite{lp19}, also arises naturally in the dual problem of D-optimal experimental design with semi-algebraic data, a standard convex optimization problem in statistics.
Numerically, problem \eqref{primal} is solved with the moment-SOS hierarchy, i.e the moment cone $\Mom_{\new{2}d}(\Spt)$ is relaxed with a hierarchy of projections of spectrahedra of increasing size.

As shown in \cite{aos19}, the Christoffel-Darboux polynomial can be used as a certificate of finite convergence of the hierarchy: the contact points of its level set at $n_d$ with $\Spt$ are the support of an optimal design. \new{An algorithmic consequence is that an optimal measure for the design problem is concentrated at the maximizers of the Christoffel-Darboux polynomial on the domain. These maximizers can be found numerically with the moment-SOS hierarchy, see \cite[Section 5.2]{aos19}.}

\new{The dual design problem has also a nice interpretation in computational geometry as an extension of the L\"owner-John ellipsopid problem to (i) semi-algebraic domains not necessarily convex and (ii) enclosing sets more general than ellipsoids.}

\new{Another interesting isssue is to study how the Christoffel-Darboux polynomial, or its maximizers, are affected when the semi-algebraic set is perturbed. It is however unclear whether our variational characterization of this polynomial can be useful for that purpose.}

The Christoffel-Darboux polynomial corresponds to a particular choice of a convex function to be minimized in the design problem, namely the logarithmic barrier function of the positive semidefinite cone. It would be \new{insightful} to study the polynomials arising in the dual design problem corresponding to other convex functions of the eigenvalues of positive semidefinite matrices \cite{lo96}.

\subsection*{Acknowledgement}
Support from the ANR-3IA Artificial and Natural Intelligence Toulouse Institute is gratefully acknowledged. \new{This work benefited from feedback from anonymous reviewers.}

\end{document}